\documentclass[a4paper]{article}
\usepackage{amsmath,amsthm,amssymb,amscd}

\numberwithin{equation}{section}

{\theoremstyle{definition}\newtheorem{definition}{Definition}[section]
\newtheorem{notation}[definition]{Notation}
\newtheorem{terminology}[definition]{Terminology}
\newtheorem{remark}[definition]{Remark}
}
\newtheorem{proposition}[definition]{Proposition}
\newtheorem{lemma}[definition]{Lemma}
\newtheorem{theorem}[definition]{Theorem}

\newcommand{\Irred}{\operatorname{Irred}}
\newcommand{\aof}{A_o{(F)}}

\newcommand{\cG}{\mathbb{G}}
\newcommand{\de}{\Delta}

\newcommand{\ot}{\otimes}
\newcommand{\wB}{\widetilde {B}}
\newcommand{\cGh}{\widehat{\mathbb{G}}}
\newcommand{\recht}{\rightarrow}
\newcommand{\B}{\operatorname{B}}

\newcommand{\sde}{\delta}
\newcommand{\vphi}{\varphi}
\newcommand{\Mor}{\operatorname{Mor}}

\newcommand{\ox}{\overline{x}}

\newcommand{\eps}{\epsilon}
\newcommand{\dimq}{\operatorname{dim}_q}
\newcommand{\ttil}{\tilde{t}}
\newcommand{\om}{\omega}

\newcommand{\Tr}{\operatorname{Tr}}

\newcommand{\mult}{\operatorname{mult}}

\newcommand{\R}{\mathbb{R}}
\newcommand{\C}{\mathbb{C}}
\newcommand{\la}{\langle}
\newcommand{\ra}{\rangle}
\newcommand{\cB}{\mathcal{B}}

\newcommand{\cst}{C$^*$}
\newcommand{\Vtil}{\widetilde{V}}

\newcommand{\deh}{\hat{\Delta}}
\newcommand{\N}{\mathbb{N}}

\newcommand{\al}{\alpha}
\newcommand{\be}{\beta}
\newcommand{\GL}{\operatorname{GL}}

\newcommand{\cK}{\mathcal{K}}

\newcommand{\SU}{\operatorname{SU}}

\newcommand{\M}{\operatorname{M}}
\newcommand{\U}{\operatorname{U}}

\newcommand{\hon}{H^{\infty}(\widehat{\mathbb{G}},\mu)}
\newcommand{\plon}{\ell^{\infty}(\widehat{\mathbb{G}})}
\newcommand{\iz}{\ensuremath{\Phi}}

\newcommand{\rP}{\mathcal{P}}

\newcommand{\id}{\mathalpha{\text{\rm id}}}
\newcommand{\End}{\operatorname{End}}
\newcommand{\fancyot}{\mathbin{\text{\footnotesize\textcircled{\tiny \sf T}}}}
\newcommand{\Vb}{\mathbb{V}}

\newcommand{\cX}{\mathbb{X}}
\newcommand{\omu}{\overline{\mu}}
\newcommand{\slim}[1]{s^*{\text{\rm -}}\lim_{#1}}
\newcommand{\epsh}{\widehat{\varepsilon}}
\newcommand{\sgn}{\operatorname{sgn}}
\newcommand{\cGtil}{\widetilde{\cG}}
\newcommand{\Htil}{\widetilde{H}}
\newcommand{\Qtil}{\widetilde{Q}}

\newcommand{\Util}{\widetilde{U}}
\newcommand{\ptil}{\widetilde{p}}
\newcommand{\psitil}{\widetilde{\psi}}
\newcommand{\Sone}{\raisebox{-0.2ex}[0ex][0ex]{$\scriptstyle S^1$}\hspace{-0.2ex}\backslash}

\newcommand{\frakB}{\mathfrak{B}}
\newcommand{\bX}{\mathbb{X}}
\newcommand{\cP}{\mathcal{P}}

\begin{document}

\begin{center}
{\LARGE\bf Identification of the Poisson and Martin boundaries of \vspace{2mm}\\
    orthogonal discrete quantum groups}

\medskip

{\sc by Stefaan Vaes$^{\text{\rm (a,b)}}$ and  Nikolas Vander Vennet$^{\text{\rm (b)}}$}
\end{center}

{\footnotesize (a) Institut de Math{\'e}matiques de Jussieu; Alg{\`e}bres d'Op{\'e}rateurs; 175,
rue du Chevaleret; F--75013 Paris (France)\\
(b) Department of Mathematics; K.U.Leuven; Celestijnenlaan 200B; B--3001 Leuven (Belgium)

e-mail: vaes@math.jussieu.fr, nikolas.vandervennet@wis.kuleuven.be}
\begin{abstract}
\noindent The Poisson and Martin boundaries for invariant random walks on the dual of the orthogonal quantum groups $A_o(F)$, are identified with
higher dimensional Podle\'s spheres that we describe in terms of generators and relations. This provides the first such identification for random
walks on non-amenable discrete quantum groups.
\end{abstract}

\section*{Introduction}

Group invariant random walks on countable groups have been studied
intensively and the identification of the associated Poisson and
Martin boundaries is a natural problem. We refer to \cite{K1} for an excellent
survey. The study of random walks on discrete quantum groups was
initiated by Biane \cite{B0} who considered duals of compact groups
and obtained a theory parallel to the theory of random walks on
discrete abelian groups.

Random walks on arbitrary discrete quantum groups (i.e.\ duals of Woronowicz' compact quantum groups \cite{wor2}) and their Poisson boundaries were
introduced by Izumi in \cite{iz1}, motivated by the study of infinite product actions of compact quantum groups. In \cite{iz1}, Izumi identified the
Poisson boundary of the dual of $\SU_q(2)$ (see \cite{wor}) with the Podle\'s sphere \cite{Podles}. Later, Neshveyev and Tuset \cite{NT} associated a
Martin boundary with a random walk on a discrete quantum group and proved that the Martin boundary of the dual of $\SU_q(2)$ is still given by the
Podle\'s sphere. This generalized Biane's work \cite{B2} on the dual of $\SU(2)$.

The Poisson boundary for the dual of $\SU_q(n)$ was computed by Izumi,
Neshveyev and Tuset in \cite{INT}, but its Martin boundary remains
mysterious. Partial results on the Martin boundary for the dual of
$\SU(n)$ were obtained by Biane \cite{B1} and Collins \cite{Collins}.

The discrete quantum groups appearing in the previous paragraphs are all duals of classical groups or their $q$-deformations. A quite different class
of compact quantum groups was introduced by Van Daele and Wang \cite{VDW} and studied by Banica \cite{banica1,banica2}. In this paper, we identify
the Poisson and Martin boundary for the dual of the orthogonal compact quantum groups $A_o(F)$.

The orthogonal compact quantum groups $A_o(F)$ are quite peculiar. Given an $n$~by~$n$ matrix $F$ satisfying $F \overline{F} = \pm 1$, they are
defined as the compact quantum group generated by an $n$-dimensional unitary representation $U$ satisfying $U = F \overline{U} F^{-1}$. On the one
hand, their representation theory is similar to the one of $\SU(2)$, both having the same fusion rules. More precisely, every $A_o(F)$ is monoidally
equivalent to an $\SU_q(2)$ for a uniquely determined $q$, in the sense of \cite{BDRV}. The quantum groups $\SU_q(2)$ appear as $A_o(F)$ for $F$ a
$2$~by~$2$ matrix. On the other hand, once $F$ is at least $3$~by~$3$, the dimensions of the irreducible representations of $A_o(F)$ start growing
exponentially, yielding a very different operator algebraic behaviour. In a sense, the operator algebras associated with $A_o(F)$ for $F$ at least
$3$~by~$3$, share several properties with the free group \cst- and von Neumann algebras, see \cite{Vbo}. In particular, the dual discrete quantum
group becomes non-amenable. As such, for the first time, Poisson and Martin boundaries of a non-amenable quantum group are identified.

In Theorem \ref{thm.main}, the Poisson boundary for the dual of $A_o(F)$ is identified with a kind of \lq higher dimensional Podle\'s sphere\rq, that
we describe in terms of generators and relations. In Theorem \ref{thm.martin}, it is shown that the Martin boundary for the dual of $A_o(F)$ can be
identified with the \cst-counterpart of these higher dimensional Podle\'s spheres.

Our method to obtain the Poisson boundary for the dual of $A_o(F)$ goes as follows. We exploit the monoidal equivalence of $A_o(F)$ and $\SU_q(2)$ in
order to reduce the identification problem to a purely $\SU_q(2)$-problem. The latter is solved invoking Izumi's computation for the Poisson boundary
of the dual of $\SU_q(2)$ (see \cite{iz1}, or the alternative approach in \cite{INT}). The Martin boundary is obtained by using a result of
\cite{Vbo}, allowing to deduce the Martin boundary from the Poisson boundary. Altogether, our proofs depend on the known computation for the Poisson
boundary of the dual of $\SU_q(2)$, but do provide an alternative method to identify the Martin boundary for the dual of $\SU_q(2)$ (as was done in
\cite{NT}).

The method of this paper to identify the Poisson boundary for the dual of $A_o(F)$ by using the notion of monoidal equivalence suggests that there is
a general way to describe the behaviour of Poisson boundaries when passing to monoidally equivalent quantum groups. This will be the subject of a
forthcoming paper of the second author.

\section{Preliminaries}

Consider a subset $S$ of a \cst-algebra. We denote by $\langle S \rangle$ the linear span of $S$ and by $[S]$ the
closed linear span of $S$.
We use the notation $\om_{\eta,\xi}(a) = \langle \eta,a \xi \rangle$ and we use inner products that are linear in the
second variable.

We use the symbol $\ot$ to denote several types of tensor products. In particular $\ot$ denotes the \emph{minimal} tensor product of \cst-algebras,
but it also denotes the tensor product of Hilbert spaces and von Neumann algebras. From the context, it will always be clear what we mean. We also
make use of the leg numbering notation in multiple tensor products: if $a \in A \ot A$, then $a_{12},a_{13},a_{23}$ denote the obvious elements in $A
\ot A \ot A$, e.g.\ $a_{12} = a \ot 1$.

\subsection*{Compact quantum groups}

We give a brief overview of the theory of compact quantum groups which
was developed by Woronowicz in \cite{wor2}. We refer to \cite{MVD} for a
survey of basic results.

\begin{definition}
A \emph{compact quantum group} $\cG$ is a pair $(C(\cG),\de)$, where
\begin{itemize}
\item $C(\cG)$ is a unital \cst-algebra; \item $\de : C(\cG) \recht C(\cG) \ot C(\cG)$ is a unital
  $^*$-homomorphism satisfying the \emph{co-associativity} relation
$$(\de \ot \id)\de = (\id \ot \de)\de \; ;$$
\item $\cG$ satisfies the \emph{left and right cancellation property}
  expressed by
$$\de(C(\cG))(1 \ot C(\cG)) \quad\text{and}\quad \de(C(\cG))(C(\cG) \ot 1)
\quad\text{are total in}\;\; C(\cG) \ot C(\cG) \; .$$
\end{itemize}
\end{definition}

\begin{remark}
The notation $C(\cG)$  suggests the analogy with the basic example given by continuous functions on a
compact group. In the quantum case however, there is no underlying space $\cG$ and $C(\cG)$ is a non-abelian
\cst-algebra.
\end{remark}

A fundamental result in the theory of compact quantum groups is the existence of a unique Haar state.
\begin{theorem}
Let $\cG$ be a compact quantum group. There exists a unique state $h$ on $C(\cG)$ which satisfies $(\id \ot h)\de(a) = h(a)1
= (h \ot \id)\de(a)$ for all $a \in C(\cG)$. The state $h$ is called the \emph{Haar state} of $\cG$.
\end{theorem}
Another crucial set of results in the framework of compact quantum groups is the Peter-Weyl representation theory.
\begin{definition}
A \emph{unitary representation} $U$ of a compact quantum group $\cG$ on a Hilbert space $H$ is a unitary element $U \in
\M(\cK(H) \ot C(\cG))$ satisfying
\begin{equation} \label{eq.rep}
(\id \ot \de)(U) = U_{12} U_{13} \; .
\end{equation}

Whenever $U^1$ and $U^2$ are unitary representations of $\cG$ on the respective Hilbert spaces $H_1$ and $H_2$, we
define
$$\Mor(U^1,U^2) := \{ T \in \B(H_2,H_1) \mid U_1(T \ot 1) = (T \ot
1)U_2 \}\; .$$ The elements of $\Mor(U^1,U^2)$ are called \emph{intertwiners}. We use the notation $\End(U) :=
\Mor(U,U)$. A unitary representation $U$ is said to be \emph{irreducible} if $\End(U) = \C1$. If $\Mor(U^1,U^2)$
 contains a unitary
operator, the
representations $U^1$ and $U^2$ are said to be \emph{unitarily equivalent}.
\end{definition}
We have the following essential result.

\begin{theorem}
Every irreducible representation of a compact quantum group is finite-dimensional. Every unitary representation is
unitarily equivalent to a direct sum of irreducibles.
\end{theorem}

Because of this theorem, we almost exclusively deal with finite-dimensional representations. By choosing an orthonormal basis of the Hilbert space $H$, a finite-dimensional unitary
representation of $\cG$ can be considered as a unitary matrix $(U_{ij})$ with entries in $C(\cG)$ and \eqref{eq.rep}
becomes
$$\de(U_{ij}) = \sum_k U_{ik} \ot U_{kj} \; .$$
The product in the \cst-algebra $C(\cG)$ yields a tensor product on
the level of unitary representations.

\begin{definition}
Let $U^1$ and $U^2$ be unitary representations of $\cG$ on the respective Hilbert spaces $H_1$ and $H_2$. We define the
tensor product
$$U^1 \fancyot U^2 := U^1_{13} U^2_{23} \in \M(\cK(H_1 \ot H_2) \ot
C(\cG)) \; .$$
\end{definition}

\begin{notation}
Let $\cG$ be a compact quantum group. We denote by $\Irred(\cG)$ the set of equivalence classes
of irreducible unitary representations.  We choose representatives $U^x$ on the Hilbert space $H_x$ for every $x \in \Irred(\cG)$. Whenever $x,y
\in \Irred(\cG)$, we use $x \ot y$ to denote the unitary representation $U^x \fancyot U^y$. The class of the trivial
unitary representation is denoted by $\varepsilon$. We define the
natural numbers $\mult(z,x \ot y)$ such that
$$x \ot y \cong \bigoplus_{z \in \Irred(\cG)} \mult(z,x \ot y) \cdot U^z \; .$$
The collection of natural numbers $\mult(z,x \ot y)$ are called the
\emph{fusion rules} of $\cG$.
\end{notation}

The set $\Irred(\cG)$ is equipped with a natural involution $x \mapsto \ox$ such that $U^{\ox}$ is the unique (up to
unitary equivalence) irreducible unitary representation satisfying
$$\Mor(x \ot \ox,\eps) \neq \{0\} \neq \Mor(\ox \ot x,\eps) \; .$$
The unitary representation $U^{\ox}$ is called the \emph{contragredient} of $U^x$.

For every $x \in \Irred(\cG)$, we take non-zero elements
${t_x}\in\Mor(x\ot\ox,\eps)$ and $s_x \in\Mor(\overline{x}\ot
x,\varepsilon)$ satisfying $(t_x^*\ot 1)(1\ot s_x)=1$. Write the
antilinear map
$$j_x : H_x\to H_{\ox}:\xi\mapsto (\xi^*\ot 1)t_x $$
and define $Q_x := j_x^*j_x$. We normalize $t_x$ in such a way that $\Tr(Q_x)=\Tr(Q_x^{-1})$. This uniquely determines $Q_x$
and fixes $t_x,s_x$ up to a number of modulus $1$. Note that $t_x^*t_x
= \Tr(Q_x)$.

\begin{definition}
For $x \in \Irred(\cG)$, the value $\Tr(Q_x)$ is called the \emph{quantum dimension} of $x$ and denoted by $\dimq(x)$.
Note that $\dimq(x) \geq \dim(x)$, with equality holding if and only if $Q_x = 1$.
\end{definition}

The irreducible representations of $\cG$ and the Haar state $h$ are connected by the \emph{orthogonality relations}.
\begin{equation}\label{eq.orthogonality}
(\id \ot h)(U^x (\xi\eta^* \ot 1) (U^y)^*) = \frac{\sde_{x,y}1}{\dimq(x)} \la \eta,Q_x\xi \ra \quad , \quad (\id \ot
h)((U^x)^*(\xi\eta^* \ot 1) U^y) = \frac{\sde_{x,y}1}{\dimq(x)} \la \eta,Q^{-1}_x\xi \ra  \; ,
\end{equation}
for $\xi\in H_x$ and $\eta\in H_y$.

\begin{notation}
Let $\cG= (C(\cG),\de)$ be a compact quantum group. We denote by $\mathcal{C}(\cG)$ the set of coefficients of finite
dimensional corepresentations of $\cG$. Hence,
\[ \mathcal{C}(\cG)=\la (\om_{\xi,\eta}\ot\id)(U^x)\mid x\in\Irred(\cG),\ \xi,\eta\in H_x
\ra \] Then, $\mathcal{C}(\cG)$ is a unital dense $^*$-subalgebra of
$C(\cG)$. Restricting $\de$ to $\mathcal{C}(\cG)$,
$\mathcal{C}(\cG)$ becomes a Hopf $^*$-algebra.\\
Also, for  $x\in \Irred(\cG)$, denote by $$\mathcal{C}(\cG)_x=\la (\om_{\xi,\eta}\ot\id)(U^x)\mid  \xi,\eta\in H_x
\ra$$
\end{notation}

\begin{definition} The \emph{reduced \cst-algebra} $C_r(\cG)$ is
  defined as the norm closure of $\mathcal{C}(\cG)$\
in the GNS-representation with respect to $h$. The \emph{universal
  \cst-algebra} $C_u(\cG)$ is defined as the enveloping \cst-algebra
of $\mathcal{C}(\cG)$. The \emph{von Neumann algebra} $L^\infty(\cG)$
is defined as the von Neumann algebra generated by $C_r(\cG)$.

A compact quantum group $\cG$ is said to be \emph{co-amenable} (and the
discrete quantum group $\widehat{\cG}$ is said to be \emph{amenable})
if the homomorphism $C_u(\cG) \recht C_r(\cG)$ is an isomorphism.
\end{definition}

Given an arbitrary compact quantum group $\cG$, we have surjective homomorphisms $C_u(\cG)\recht C(\cG)\recht C_r(\cG)$, but we are only interested
in $C_r(\cG)$ and $C_u(\cG)$. Note that if $\cG$ is the dual of a
discrete group $\Gamma$, we have $C_r(\cG) = C^*_r(\Gamma)$ and
$C_u(\cG) = C^*(\Gamma)$.

\begin{proposition} The Haar state $h$ is a KMS-state on both $C_r(\cG)$\
  and $C_u(\cG)$ and the modular group is determined by $$(\id\ot
\sigma_t^h)(U^x)=(Q^{it}_x\ot 1)U^x(Q^{it}_x\ot 1)$$ for every $x\in \Irred(\cG)$.
\end{proposition}

\subsection*{Discrete quantum groups and duality}

A discrete quantum group is defined as the dual of a compact quantum group by putting together all irreducible
representations.

\begin{definition}
Let $\cG$ be a compact quantum group. We define the dual (discrete) quantum group $\cGh$ as follows.
\begin{equation*}
c_0(\cGh) = \bigoplus_{x \in \Irred(\cG)} \B(H_x) \; , \qquad
\ell^\infty(\cGh) = \prod_{x \in \Irred(\cG)} \B(H_x) \; .
\end{equation*}
We denote the minimal central projections of $\ell^\infty(\cGh)$ by $p_x$, $x \in \Irred(\cG)$. We have a natural
unitary $\Vb \in \M(c_0(\cGh) \ot C(\cG))$ given by
\begin{equation}\label{UM}
\Vb = \bigoplus_{x \in \Irred(\cG)} U^x \; .
\end{equation}
This unitary $\Vb$ implements the duality between $\cG$ and $\cGh$. We have a natural comultiplication
$$\deh : \ell^\infty(\cGh) \recht
\ell^\infty(\cGh) \ot \ell^\infty(\cGh) : (\deh \ot \id)(\Vb) = \Vb_{13} \Vb_{23} \; .$$
\end{definition}

One can deduce from this the following equivalent way to define the product structure on $\ell^\infty(\cGh)$.
$$\deh(a) S =
S a \quad\text{for all}\;\; a \in \ell^\infty(\cGh),\  S \in \Mor(y \ot z,x) \; .$$  The notation introduced above is
aimed to suggest the basic example where $\cG$ is the dual of a discrete group $\Gamma$, given by $C(\cG) =
C^*(\Gamma)$ and $\de(\lambda_x) = \lambda_x \ot \lambda_x$ for all $x \in \Gamma$. The map $x \mapsto \lambda_x$
yields an identification of $\Gamma$ and $\Irred(\cG)$ and then, $\ell^\infty(\cGh) = \ell^\infty(\Gamma)$.

\begin{remark}
It is of course possible to give an intrinsic definition of a discrete quantum group. This was already implicitly clear in Woronowicz' work and was explicitly done in \cite{ER,VD}. For our
purposes, it is most convenient to take the compact quantum group as a starting point: indeed, all interesting examples
of concrete discrete quantum groups are defined as the dual of certain compact quantum groups.
\end{remark}

The discrete quantum group $\ell^\infty(\cGh)$ comes equipped with a natural modular structure.

\begin{notation} \label{not.states}
 We have canonically defined states $\vphi_x$ and $\psi_x$ on $\B(H_x)$
related to \eqref{eq.orthogonality} as follows.
\begin{align*}
\psi_x(A) &= \frac{1}{\dimq(x)}t_x^* (A \ot 1) t_x = \frac{\Tr(Q_xA)}{\Tr(Q_x)} = (\id \ot h)(U^x (A \ot 1) (U^x)^*) \quad\text{and} \\
\vphi_x(A) &=\frac{1}{\dimq(x)} t_{\ox}^* (1 \ot A) t_{\ox} = \frac{\Tr(Q_x^{-1}A)}{\Tr(Q_x^{-1})} = (\id \ot
h)((U^x)^*(A \ot 1) U^x) \; ,
\end{align*}
for all $A \in \B(H_x)$.
\end{notation}
\begin{remark}
The states $\vphi_x$ and $\psi_x$ are significant, since they provide a formula for the invariant weights on
$\ell^\infty(\cGh)$. The left invariant weight is given by  $\sum_{x \in \Irred(\cG)} \dimq(x)^2 \psi_x $, and the
right invariant weight is
 given by $\sum_{x \in \Irred(\cG)} \dimq(x)^2 \vphi_x$.
 \end{remark}

\subsection*{Examples: the universal orthogonal compact quantum groups}

This class of compact quantum groups was introduced by Wang and Van
Daele in \cite{VDW} and studied by Banica in \cite{banica1}.

\begin{definition}
Let $F \in \GL(n,\C)$ satisfying $F \overline{F} = \pm 1$. We define the compact quantum group $\cG = A_o(F)$ as
follows.
\begin{itemize}
\item $C(\cG)$ is the universal \cst-algebra with generators
  $(U_{ij})$ and relations making $U = (U_{ij})$ a unitary element of
  $\M_n(\C) \ot C(\cG)$ and $U = F \overline{U} F^{-1}$, where $(\overline{U})_{ij} = (U_{ij})^*$.
\item $\de(U_{ij}) = \sum_k U_{ik} \ot U_{kj}$.
\end{itemize}
\end{definition}

In these examples, the unitary matrix $U$ is a representation, called the \emph{fundamental representation}. The
definition of $\cG=A_o(F)$ makes sense without the requirement $F \overline{F} = \pm 1$, but the fundamental
representation is irreducible if and only if $F \overline{F} \in \R
1$. We then normalize such that $F \overline{F} = \pm 1$.

\begin{remark}
It is easy to classify the quantum groups $A_o(F)$. For $F_1,F_2 \in \GL(n,\C)$ with $F_i \overline{F}_i = \pm 1$, we
write $F_1 \sim F_2$ if there exists a unitary matrix $v$ such that $F_1 = v F_2 v^t$, where $v^t$ is the transpose of
$v$. Then, $A_o(F_1) \cong A_o(F_2)$ if and only if $F_1 \sim F_2$. It follows that the $A_o(F)$ are classified up to
isomorphism by $n$, the sign $F \overline{F}$ and the eigenvalue list of $F^*F$ (see e.g.\ Section 5 of \cite{BDRV}
where an explicit fundamental domain for the relation $\sim$ is described).

If $F \in \GL(2,\C)$, we get up to equivalence, the matrices
\begin{equation}\label{eq.fq}
F_q = \begin{pmatrix} 0 & |q|^{1/2} \\ - (\sgn q) |q|^{-1/2} & 0
\end{pmatrix}
\end{equation}
for $q \in \, [-1,1]$, $q \neq 0$, with corresponding quantum groups $A_o(F_q) \cong \SU_{q}(2)$, see \cite{wor}.
\end{remark}

The following result has been proved by Banica \cite{banica1}. It
tells us that the compact quantum groups $A_o(F)$ have the same fusion
rules as the group $\SU(2)$.

\begin{theorem}
Let $F \in \GL(n,\C)$ and $F \overline{F} = \pm 1$. Let $\cG = A_o(F)$. Then  $\Irred(\cG)$ can be identified with $\N$ in
such a way that
$$x \ot y \cong |x-y| \oplus (|x-y|+2) \oplus \cdots \oplus (x+y) \;
,$$ for all $x,y \in \N$.
\end{theorem}

\subsection*{Actions of quantum groups and spectral subspaces}

\begin{definition}\label{actie} Let $B$ be a unital \cst-algebra. A (right) action
  of $\cG$ on $B$ is a unital $^*$-homomorphism
$\sde : B \recht B \ot C(\cG)$ satisfying
\[(\sde\ot \id)\sde = (\id \ot \de)\sde \quad\text{and}\quad [\sde(B)(1\ot C(\cG))]=B\ot C(\cG)\; .\] The
action $\sde$ is said to be ergodic if the fixed point algebra $B^{\sde}:=\{x\in B\mid  \sde (x)=x\ot 1\}$ equals $\C
1$. In that case, $B$ admits a unique invariant state $\om$ given by $\om(b) 1 = (\id \ot h)\sde(b)$.
\end{definition}

\begin{definition} \label{def.spectral}
Let $\delta : B \recht B \ot C(\cG)$ be an action of the compact
quantum group $\cG$ on the unital \cst-algebra $B$.
\begin{itemize}
\item For every $x \in \Irred(\cG)$, the spectral subspace $B_x$ is
  defined as the linear subspace $B_x \subset B$ given by
\begin{equation*}
B_x:= \{b\in B\mid\ \sde(b)\in B\ot \mathcal{C}(\cG)_x \}
\end{equation*}
Note that $\delta : B_x \recht B_x \ot \mathcal{C}(\cG)_x$ and that
$B_x^* = B_{\ox}$.
\item We define $\cB$ as the linear span of the $B_x$, $x \in
  \Irred(\cG)$. Then, $\cB$ is a dense $^*$-subalgebra of $B$.
\item The action $\delta$ is said to be \emph{universal} if $B$ is the
  universal enveloping \cst-algebra of $\cB$. It is said to be
  \emph{reduced}  if the
  conditional expectation $(\id\ot h)\sde$ of $B$ on $B^\sde$ is faithful.
\item If $\delta$ is ergodic, $B_x$ is finite dimensional and its
  dimension is of the form $\dim H_x \cdot \mult(x,\delta)$, where
  $\mult(x,\delta)$ is called the multiplicity of $x$ in $\delta$.
\end{itemize}
\end{definition}

\begin{remark}\label{amenable}
If $\cG$ is co-amenable, $C(\cG)$ has a bounded co-unit and a faithful
Haar state. Hence, any action of $\cG$ is both universal and reduced.
\end{remark}

Actions on von Neumann algebras are defined as follows.

\begin{definition} A right action of a compact (resp.\ discrete)
  quantum group $\cG$ (resp. $\widehat{\cG}$) on a von Neumann algebra $N$ is an injective
 normal unital $^*$-homomorphism
$$\sde:N\recht N\ot L^{\infty}(\cG)\qquad \text{resp.}\quad \sde:N\recht N\ot \ell^{\infty}(\widehat{\cG})$$
satisfying $(\sde\ot \id)\sde = (\id \ot \de)\sde $, resp.\ $(\delta
\ot \id)\delta = (\id \ot \deh)\delta$.
\end{definition}

\section{The Poisson boundary of a discrete quantum group}

We give a brief survey of Izumi's theory of Poisson boundaries for
discrete quantum groups and his
computation of the Poisson boundary for the dual of $\SU_q(2)$, see
\cite{iz1}. We state as well some of the results of \cite{INT} that
are useful in the rest of the paper.

\subsection*{General results}

Fix a discrete quantum group $\widehat{\cG}$.

\begin{notation}
For every normal state $\phi\in \plon$, we define the convolution operator
 $$P_{\phi}:\plon\recht\plon:P_{\phi}(a)=(\id \ot
 \phi)\widehat{\Delta}(a) \; .$$
\end{notation}

We are only interested in special states $\phi\in \plon$, motivated by
the following straightforward proposition. For every probability measure $\mu$ on
$\Irred(\cG)$, we set
\begin{equation*}
\psi_{\mu}=\sum_{x\in \Irred(\cG)}\mu(x)\psi_x \quad\text{and}\quad P_\mu
:= P_{\psi_{\mu}} \; .
\end{equation*}
Recall that the states $\psi_x$ are defined in \ref{not.states}. Note that we have a convolution product $\mu * \nu$ on the
measures on $\Irred(G)$, such that $\psi_{\mu * \nu} = (\psi_\mu \ot \psi_\nu)\deh$.

\begin{proposition}
Let $\phi$ be a normal state on $\plon$. Then the following conditions are equivalent. \begin{itemize}
\item $\phi$ has the form $\psi_{\mu}$ from some probability measure $\mu$ on
$\Irred(\cG)$.
\item The Markov operator $P_\phi$ preserves the center of $\ell^\infty(\cGh)$.
\item $\phi$ is invariant under the adjoint action $\al_\cG :
  \ell^\infty(\cGh) \recht \ell^\infty(\cGh) \ot L^\infty(\cG) : a \mapsto
  \mathbb{V}(a\ot 1)\mathbb{V}^*$.
\end{itemize}
\end{proposition}

\begin{definition}[\cite{iz1}, Section 2.5]
Let $\mu$ be a probability measure on $\Irred(\cG)$. Set
 $$H^{\infty}(\widehat{\cG},\mu)=\{a \in \plon\mid\  P_{\mu}(a)=a\}
 \; .$$
Equipped with the product defined by
$$a \cdot b := \slim{n \recht \infty} \frac{1}{n} \sum_{k=1}^n
P_\mu^k(ab) \; ,$$
the space $H^{\infty}(\widehat{\cG},\mu)$ becomes a von Neumann
algebra that we call the \emph{Poisson boundary} of $\cGh$ with
respect to $\mu$.
\end{definition}

\begin{terminology}
A probability measure
$\mu$ on $\Irred({\cG})$ is called \emph{generating} if there exists, for every $x\in \Irred({\cG})$, an $n\geq 1$ such that $\mu^{* n}(x)\not=0$.
\end{terminology}

In all the results below, the probability measure $\mu$ is assumed to
be generating. In that case, it is not hard to concretely realize
$H^\infty(\cGh,\mu)$ as a von Neumann algebra. Indeed,
$$\pi_\infty : H^\infty(\cGh,\mu) \recht \bigotimes_{n=0}^\infty
(\ell^\infty(\cGh),\psi_\mu) : \pi_\infty(a) = \slim{n \recht \infty}
\pi_n(a) \quad\text{where}\;\; \pi_n := \deh^{(n)} : \ell^\infty(\cGh)
\recht \bigotimes_{k=0}^{n-1} \ell^\infty(\cGh) \; .$$
Moreover, the product becomes $a \cdot b = \slim{n \recht \infty} P_\mu^n(ab)$.

The restriction of the co-unit $\epsh$ yields a state on $\hon$,
called the \emph{harmonic state}. It is
clear that this state is faithful when $\mu$ is generating.

\begin{definition} \label{def.actions}
Let $\mu$ be a generating measure on $\Irred(\cG)$. The Poisson
boundary $H^\infty(\cGh,\mu)$ comes equipped with two natural actions,
one of $\cG$ and one of $\cGh$:
\begin{align*}
\al_{\cG} & : \hon \recht \hon \ot L^\infty(\cG) : \al_{\cG}(a) =
\mathbb{V}(a\ot 1)\mathbb{V}^* \; ,  \\
\al_{\cGh} & : \hon \recht \ell^\infty(\cGh) \ot \hon : \al_{\cGh}(a) =
  \deh(a) \; .
\end{align*}
\end{definition}
Note that $\al_{\cG}$ is the restriction of the adjoint action of
$\cG$ on $\ell^\infty(\cGh)$, while $\al_{\cGh}$ is nothing else than
the restriction of the comultiplication. The actions $\al_{\cG}$ and
$\al_{\cGh}$ are well defined because of the following equivariance
formulae:
$$(\id\ot P_{\mu})(\widehat{\Delta}(a))=\widehat{\Delta}(P_{\mu}(a))
\quad\text{and}\quad (P_{\mu}\ot \id)(\al_\cG(a))=\al_\cG(P_{\mu}(a)) \; .$$
When $\widehat{\cG}$ is a discrete group, the action $\al_{\cG}$ is the trivial action
on $\plon$.
In general, the fixed point algebra of $\al_{\cG}$
is precisely the algebra of \emph{central harmonic elements}
$Z(\plon)\cap \hon$. Since the Markov operator $P_{\mu}$ preserves the
center $Z(\plon)$, the commutative von Neumann algebra $Z(\plon)\cap
\hon$ with state $\epsh$, is exactly the Poisson boundary for the random
walk on $\Irred(\cG)$ with transition probabilities
$p(x,y)$ and $n$-step transition probabilities $p_n(x,y)$
given by
\begin{equation}\label{eq.nstep}
p_x p(x,y) = p_x P_\mu(p_y) \; , \quad p_x p_n(x,y) = p_x P_\mu^n(p_y)
\; .
\end{equation}
Note that $p_n(e,y) = \mu^{* n}(y) = \psi_\mu^{*
n}(p_y)$.

So, the action $\al_{\cG}$ is \emph{ergodic} if and
only if there are no non-trivial central harmonic elements. This
occurs if the fusion algebra of $\cG$ is commutative. Much more can be
said in that case. We record the following results for future use.

\begin{proposition}\label{prop.abelian}
Suppose that the fusion algebra of $\cG$ is commutative
(i.e.\ $\mult(y\ot z,x)=\mult(z\ot y,x)$ for every $x,y,z \in
\Irred(\cG)$) and let $\mu$ be a generating probability measure on $\Irred(\cG)$.
\begin{itemize}
\item (Cor.\ 3.5 in \cite{Hai} and Cor.\ 3.2 in \cite{INT}) There are
  no non-trivial central harmonic elements, i.e.\ $$Z(\plon)\cap
  \hon=\mathbb{C}1 \; .$$
\item (Prop.\ 1.1 in \cite{INT}) The Poisson boundary does not depend
  on the choice of generating measure:
$$\hon = \{a \in \ell^\infty(\cGh) \mid P_x(a) = a \quad\text{for
  all}\;\; x \in \Irred(\cG) \} \; .$$
\item (Cor.\ 3.5 in \cite{INT}) Using the notation of \ref{def.spectral}, we have
$$\mult(x,\alpha_{\cG})\leq \sup\{\mult(\overline{y}\ot y,x )\ \vert\ y\in
\Irred(\cG)\} \; ,$$
for all $x \in \Irred(\cG)$.
\end{itemize}
\end{proposition}

\subsection*{The Poisson boundary of $\widehat{SU_q(2)}$}

Let $\cG=\SU_q(2)$ for $-1 < q < 0$ or $0 < q < 1$. In \cite{iz1},
Izumi identified the Poisson boundary of $\cGh$ with the Podle\'s
sphere. Since this result is needed in our identification of the
Poisson boundary of the dual of $A_o(F)$, we briefly recall it
here. We also refer to \cite{INT} for an easier approach to the
computation of the Poisson boundary of the dual of $\SU_q(2)$.

We first record the following general result.

\begin{proposition}[Lemma 3.8 in \cite{iz1}] \label{prop.izop}
Let $\cG$ be a compact quantum group. Defining
\begin{equation*}
{\iz}:L^{\infty}(\cG)\recht
\ell^{\infty}(\widehat{\cG}):\iz(a)=(\id\ot {h})(\mathbb{V}^*(1\ot
a)\mathbb{V}) \; ,
\end{equation*}
the image of $\iz$ is contained in $\hon$ for any probability measure
$\mu$ on $\Irred(\cG)$.
\end{proposition}

For the rest of this section, fix $q \in \, ]-1,1[\, , q \neq 0$ and set $\cG = \SU_q(2)$ (see \cite{wor}). Recall that $C(SU_q(2))$ is generated by
the coefficients of the unitary operator $U = \bigl(\begin{smallmatrix} \alpha & -q\gamma^* \\ \gamma &
  \alpha^* \end{smallmatrix}\bigr)$ and that $\Delta$ is defined such
that $U$ is a unitary representation of $\SU_q(2)$, called the
fundamental representation.

Consider the natural homomorphism
$$\pi_{S^1}:\SU_q(2)\rightarrow C(S^1):\pi_{S^1}(\alpha)=z\quad \mbox{and}\quad \pi_{S^1}(\gamma)=0$$
So, $\Delta_{S^1}\pi_{S^1}=(\pi_{S^1}\otimes \pi_{S^1})\Delta$, and we can consider $\pi_{S^1}$ as
 an embedding of the circle $S^1$ into $\SU_q(2)$.\\
Define
 $$C(S^1\backslash\SU_q(2)):=C(\SU_q(2))^{S^1}=\{a\in \SU_q(2) \mid
 (\pi_{S^1}\ot 1){\Delta}(a)=1\ot a\}$$
 This homogeneous space is just the fixed point algebra of the action $$\Delta_{S^1}:=(\pi_{S^1}\ot \id)\Delta :C(SU_q(2))\recht C(S^1)\ot C(SU_q(2))$$
  and is called the Podle\'s sphere \cite{Podles}.

Observe that the
 restriction of the comultiplication yields a right action
 $\beta_{\cG}$ of
 $\cG$ on $L^\infty(S^1\backslash \SU_q(2))$. Moreover, the restriction
 of the adjoint action
$$\beta_{\cGh} : L^\infty(\cG) \recht \ell^\infty(\cGh) \ot
L^\infty(\cG) : \beta_{\cGh}(a) = \mathbb{V}^*(1\ot a)\mathbb{V}$$
yields a left action of $\cGh$ on $C(S^1\backslash \SU_q(2))$.

\begin{theorem}[Th.\ 5.10 in \cite{iz1}] \label{psuq}
Let $\cG = \SU_q(2)$ ($q \neq \pm 1$) and let $\mu$ be a generating probability measure on $\Irred(\cG)$. The
restriction of the completely positive map ${\iz}$ (see Proposition
\ref{prop.izop}) to $L^{\infty}(S^1\backslash \cG)$ is a
$^*$-isomorphism $\iz_0$
between the Podle\'s sphere and the Poisson boundary of $\cGh$.

The $^*$-isomorphism $\iz_0$ intertwines the actions
$\be_{\cG},\be_{\cGh}$ on $L^\infty(S^1 \backslash \cG)$ with the
actions $\al_{\cG},\al_{\cGh}$ on $\hon$ defined
in \ref{def.actions}. Finally, $\iz_0$ intertwines the harmonic state $\epsh$
with the restriction of the Haar state.
\end{theorem}

\section{The Martin boundary of a discrete quantum group}

The Martin boundary and the Martin compactification of a discrete
quantum group have been defined by Neshveyev and Tuset in
\cite{NT}.
Fix a discrete quantum group $\cGh$ and
a probability measure $\mu$ on $\Irred(\cG)$. We have an
associated Markov operator $P_\mu$ and a classical random walk on $\Irred(\cG)$ with $n$-step transition probabilities
given by \eqref{eq.nstep}.

\begin{definition}
The probability measure $\mu$ on $\Irred(\cG)$ is said to be \emph{transient}
if $\sum_{n=0}^\infty p_n(x,y) < \infty$ for all $x,y \in \Irred(\cG)$.
\end{definition}

We suppose throughout that $\mu$ is a
generating measure and that $\mu$ is transient.

Denote by $c_c(\cGh) \subset c_0(\cGh)$ the algebraic direct sum of
the algebras $\B(H_x)$. We define, for $a \in c_c(\cGh)$,
$$G_\mu(a) = \sum_{n=0}^\infty P_\mu^n(a) \; .$$
Observe that usually $G_\mu(a)$ is unbounded, but it makes sense in
the multiplier algebra of $c_c(\cGh)$, i.e.\ $G_\mu(a)p_x \in \B(H_x)$
makes sense for every $x \in \Irred(\cG)$ because $\mu$ is transient. Moreover, $G_\mu(p_\eps)$ is strictly positive and
central. This allows to define the Martin kernel as follows.

Whenever $\mu$ is a measure on $\Irred(\cG)$, we use the notation
$\omu$ to denote the measure given by $\omu(x) = \mu(\ox)$.

\begin{definition}[Defs.\ 3.1 and 3.2 in \cite{NT}]
Define
$$K_\mu : c_c(\cGh) \recht \ell^\infty(\cGh) : K_\mu(a) = G_\mu(a)
G_\mu(p_\eps)^{-1} \; .$$
Define the \emph{Martin compactification} $\widetilde{A}_{\mu}$ as the
\cst-subalgebra of $\ell^\infty(\cGh)$ generated by
$K_{\omu}(c_c(\cGh))$ and $c_0(\cGh)$. Define the \emph{Martin
  boundary} $A_{\mu}$ as the quotient $\widetilde{A}_{\mu} / c_0(\cGh)$.
\end{definition}

By Theorem 3.5 in \cite{NT}, the adjoint action $\al_\cG$ and the
comultiplication $\deh$ define, by restriction and passage to the
quotient, the following actions on the Martin boundary.
\begin{equation}\label{eq.actionsMartin}
\gamma_\cG : A_\mu \recht A_\mu \ot C(\cG) \quad\text{and}\quad
\gamma_{\cGh} : A_\mu \recht \M(c_0(\cGh) \ot A_\mu) \; .
\end{equation}

\section{Monoidal equivalence}

A crucial tool in the computation of the Poisson boundary for $\cGh$ when $\cG = A_o(F)$ is the \emph{monoidal equivalence} of $A_o(F)$ and
$\SU_q(2)$ for the appropriate value of $q$. This notion was introduced in \cite{BDRV} and is reviewed in this section.

\begin{definition}[Def.\ 3.1 in \cite{BDRV}] \label{def.moneq}
Two compact quantum groups $\cG_1=(C(\cG_1),\de)$ and $\cG_2=(C(\cG_2),\de_2)$ are said to be \emph{monoidally
equivalent} if there exists a bijection $\vphi:\cGh_1\to\cGh_2$ satisfying $\vphi(\varepsilon) = \varepsilon$, together
with linear isomorphisms
\[\vphi:\Mor(x_1 \ot \cdots \ot x_r ,y_1\ot\cdots\ot
y_k)\to\Mor(\vphi(x_1) \ot \cdots \ot \vphi(x_r),\vphi(y_1)\ot\cdots \ot \vphi(y_k))\] satisfying the following
conditions:
\begin{equation*}
\begin{alignedat}{2}
\vphi(1) &= 1 & \qquad
\vphi(S \ot T) &= \vphi(S) \ot \vphi(T) \\
\vphi(S^*) &= \vphi(S)^* & \qquad \vphi(S T) &=\vphi(S) \vphi(T)
\end{alignedat}
\end{equation*}
whenever the formulas make sense. In the first formula, we consider $1 \in \Mor(x,x) = \Mor( x \ot \varepsilon,x) =
\Mor(\varepsilon \ot x,x)$. Such a collection of maps $\vphi$ is called a \emph{monoidal equivalence} between $\cG_1$
and $\cG_2$.
\end{definition}

By Theorem 3.9 and Proposition 3.13  of \cite{BDRV}, we have the following.

\begin{theorem}\label{construction}
Let $\vphi$ be a monoidal equivalence between compact quantum groups $\cG_1$ and $\cG_2$.
\begin{itemize}
\item There exist a unique unital $^*$-algebra $\cB$ equipped with a faithful state $\om$ and unitary elements $X^x \in
\B(H_x,H_{\vphi(x)}) \ot \cB$ for all $x \in \cGh_1$, satisfying
\begin{enumerate}
\item $X^y_{13} X^z_{23} (S \ot 1) = (\vphi(S) \ot 1)X^x
  \quad\text{for all}\quad S \in \Mor(y\ot z,x) \; ,$
\item the matrix coefficients of the $X^x$ form a linear basis of
  $\cB$,
\item $(\id \ot \om)(X^x) = 0 \quad\text{if}\quad x \neq \eps$.
\end{enumerate}
\item There exists  unique commuting ergodic actions $\sde_1 : \cB \recht \cB \ot C(\cG_1)$ and $\sde_2 : \cB \recht C(\cG_2)\ot
\cB$ satisfying
$$(\id \ot \sde_1)(X^x) = X^x_{12} U^x_{13}\qquad \mbox{and}\qquad (\id \ot \sde_2)(X^x) = U^{\varphi(x)}_{12} X^x_{13}$$
for all $x \in \cGh$.
\item The state $\om$ is invariant under $\sde_1$ and $\sde_2$. Denoting
  by $B_r$ the \cst-algebra generated by $\cB$ in the
  GNS-representation associated with $\om$ and denoting by $B_u$ the
  universal enveloping \cst-algebra of $\cB$, the actions $\sde_1,\sde_2$
  admit unique extensions to actions on $B_r$ and $B_u$.
\end{itemize}
\end{theorem}

Note that in the case $\cG=\cG_1=\cG_2$ and $\vphi$ the identity map, we have $\cB=\mathcal{C}(\cG)$ and $X^x=U^x$ for every $x\in \Irred(\cG)$. The
following unitary operator generalizes \eqref{UM}.
\begin{equation}\label{eq.X}
\mathbb{X}:=\bigoplus_{x\in\Irred(\cG)}X^x \quad\text{where}\quad
\mathbb{X} \in \prod_{x \in \Irred(\cG)} \bigl(\B(H_x,H_{\vphi(x)})
\ot B \bigr) \; .
\end{equation}

\begin{proposition}
The invariant state $\om$ is a KMS state on $B_r$ and $B_u$ and its modular group is determined by
\begin{equation}\label{KMSS}
(\id\ot \sigma^{\om}_t)(X^x)=(Q_{\vphi(x)}^{it}\ot 1)X^x(Q_x^{it}\ot 1)
\end{equation}
for every $x\in \Irred(\cG_1)$.
\end{proposition}

\begin{remark}
Define $B_x:=\la (\om_{\xi,\eta}\ot \id)(X^x)\mid\ \xi\in H_{\varphi(x)}, \eta\in H_x
\ra \; .$ Then, as a vector space $$\mathcal{B}=\bigoplus_{x\in \Irred(\cG)}B_x\; .$$ Moreover, the $B_x$ are exactly the spectral subspaces of
$\sde_1$ and $\sde_2$, while $\mathcal{B}$ is exactly the dense $^*$-algebra given by Definition \ref{def.spectral}.
\end{remark}

The orthogonality relations \eqref{eq.orthogonality} generalize and take the following form.
\begin{equation}\label{eq.orthogonalitymon}
\begin{split}
(\id \ot \om)(X^x (\xi_1\eta_1^* \ot 1) (X^y)^*) &= \frac{\sde_{x,y}1}{\dimq(x)} \la \eta_1,Q_x\xi_1 \ra \; , \\ (\id \ot \om)((X^x)^*((\xi_2\eta_2^*
\ot 1) X^y) &= \frac{\sde_{x,y}1}{\dimq(x)}\la \eta_2,Q_{\varphi(x)}^{-1}\xi_2 \ra \; ,
\end{split}
\end{equation}
for $\xi_1\in H_x$, $\eta_1\in H_y$, $\xi_2\in H_{\varphi(x)}$ and  $\eta_2\in H_{\varphi(y)}$.

We turn to the case where $\cG_1=A_o(F_1)$ and $\cG_2=A_o(F_2)$, which will be needed in this article.

\begin{theorem}[Thms.\ 5.3 and 5.4 in \cite{BDRV}] \label{constructionk}
 Let $F_1\in M_{n_1}(\mathbb{C})$ and $F_2\in M_{n_2}(\mathbb{C})$ such that $F_1\overline{F}_1=\pm 1$\
 and $F_2\overline{F}_2=\pm 1$.
\begin{itemize}
\item  The compact quantum groups $A_o(F_1)$ and $A_o(F_2)$
  are monoidally equivalent iff $F_1\overline{F}_1$ and $F_2\overline{F}_2$ have the same sign and
  $\Tr(F_1^*F_1)=\Tr(F_2^*F_2)$.
\item Assume that $A_o(F_1)$ and $A_o(F_2)$
  are monoidally equivalent.
  Denote by $A_o(F_1,F_2)$ the universal unital \cst-algebra
  generated by the coefficients of
$$Y \in M_{n_2,n_1}(\C) \ot A_o(F_1,F_2) \quad\text{with relations}\quad
Y \;\;\text{unitary}\quad\text{and}\quad Y = (F_2 \ot 1)\overline{Y} (F_1^{-1} \ot 1) \; .$$ Then, $A_o(F_1,F_2) \neq 0$ and there exists a unique
pair of commuting universal ergodic actions, $\sde_1$ of $A_o(F_1)$ and $\sde_2$ of $A_o(F_2)$, such that $$(\id \ot \sde_1)(Y) = Y_{12} (U_1)_{13}
\quad\text{and}\quad (\id \ot \sde_2)(Y) = (U_2)_{12} Y_{13} \; .$$ Here, $U_i$ denotes the fundamental corepresentation of $A_o(F_i)$.
\item $(A_o(F_1,F_2),\sde_1,\sde_2)$ is isomorphic with the \cst-algebra $B_u$ and the actions thereon given by Theorem \ref{construction}
\end{itemize}
\end{theorem}

\section{Poisson boundary of the dual of $\aof$}

Fix $n \geq 2$ and a matrix $F \in \GL(n,\C)$ satisfying $F \overline{F} = \pm 1$. Set $\cG = A_o(F)$ which remains fixed throughout this section. We
assume that $\cG \neq \SU_{\pm 1}(2)$. We identify the Poisson boundary $H^\infty(\cGh,\mu)$ for a generating measure $\mu$ on $\Irred(\cG)$.

Take the unique $q \in \, ]-1,1[$ such that $F \overline{F} = - \sgn q$ and $\Tr(F^*F) = |q + 1/q|$. By Theorem \ref{constructionk}, $\cG$ is
monoidally equivalent with $\SU_q(2)=A_o(F_q)$, where $F_q$ is given by \eqref{eq.fq}.

\begin{definition} \label{def.B}
We define the \cst-algebra $B:=A_o(F,F_q)$ and denote by $\beta_{\cG} : B \recht B \ot C(\cG)$ the ergodic action of $\cG$ with invariant state
$\om$, given by Theorem \ref{constructionk}. So, $B$ is generated by the entries of a unitary $2$~by~$n$ matrix $Y$ satisfying $Y = F_q \overline{Y}
F^{-1}$. Define the action $\rho$ of $S^1$ on $B$ such that $\rho_z(Y) = \bigl(\begin{smallmatrix} z & 0 \\ 0 & \overline{z}
\end{smallmatrix}\bigr)Y$ and set
$$\Sone B:=\{x\in B \mid \rho_z(x) = x \quad\text{for all}\;\; z \in S^1 \} \; .$$
The von Neumann algebra generated by $B$ in the GNS-construction for
$\om$ is denoted by $\mathfrak{B}:=(B,\om)''$. By Remark
\ref{amenable} (and the co-amenability of $\SU_q(2)$), the state $\om$
is a faithful KMS state on $B$ and we regard $B$ as a dense subalgebra
of $\frakB$.

The generalized Izumi operator is defined as
\begin{equation}\label{eq.geniz}
\iz:{\mathfrak B}\recht \plon:\iz(a)=(\id \ot\omega)(\cX^*(1\ot a)\cX) \; ,
\end{equation}
where $\cX$ is given by \eqref{eq.X}.
\end{definition}
Note that, by Theorem \ref{constructionk}, the quantum group $\SU_q(2)$ admits a (left) ergodic action on $B$. The action $(\rho_z)_{z \in S^1}$ is
nothing else than the restriction of that ergodic action to the closed subgroup $S^1$ of $\SU_q(2)$. In this way, $\Sone B$ is a higher
dimensional counterpart of the Podle\'s sphere.

Apart from the ergodic action $\beta_{\cG}$ of $\cG$ on $B$, we also have the analogue of the adjoint action, defined as
\begin{equation}\label{eq.adjoint}
\beta_{\cGh} : \frakB \recht \ell^\infty(\cGh) \ot \frakB : \beta_{\cGh}(x) = \bX^* (1 \ot x) \bX \; ,
\end{equation}
where again, $\bX$ is given by \eqref{eq.X}. It can be checked easily
that $\beta_{\cG}$ and $\beta_{\cGh}$ leave globally invariant $\Sone
\frakB$, yielding actions on $\Sone \frakB$ that we still denote by
$\be_\cG$, $\be_{\cGh}$.

The following is the main result of the paper.
\begin{theorem} \label{thm.main}
Let $\cG = A_o(F)$ and let $q$ be as above. Let $\mu$ be a generating measure on $\Irred(\cG)$. The restriction of the generalized Izumi operator given
by \eqref{eq.geniz} yields a $^*$-isomorphism $\iz : \Sone \frakB \recht H^\infty(\cGh,\mu)$. This $^*$-isomorphism intertwines the actions
$\beta_{\cG},\beta_{\cGh}$ (see Def.\ \ref{def.B} and formula \eqref{eq.adjoint}) with the actions $\al_{\cG},\al_{\cGh}$ defined in
\ref{def.actions}.
\end{theorem}

\begin{remark}
In order to get a better understanding of the \cst-algebra $B$ and the different actions on it, one should look at the case $F = F_q$. Then, $B=
C(\SU_q(2))$ and $\Sone B$ is exactly the Podle\'s sphere. Also, the generalized Izumi operator \eqref{eq.geniz} coincides with the Izumi operator of
Proposition \ref{prop.izop}. Nevertheless, our proof of Theorem \ref{thm.main} does not provide an alternative way of identifying the Poisson
boundary for the dual of $\SU_q(2)$, because Izumi's theorem \ref{psuq} is an ingredient of our proof.
\end{remark}

The proof of Theorem \ref{thm.main} is given at the end of the section, as a combination of several preliminary results.
\begin{itemize}
\item We take a closer look at the generalized Izumi-operator $\iz$ given by \eqref{eq.geniz} and prove the equivariance of $\iz$ with respect to the actions $\be_\cG$,
$\be_{\cGh}$ and $\al_\cG$, $\al_{\cGh}$.
\item By definition, the generalized Izumi-operator $\iz$ is a normal unital completely positive mapping and we prove that $\iz$ is
multiplicative on $\Sone {\mathfrak B}$. The proof uses a technique of \cite{INT}, which allows to reduce to the case of $\cG = \SU_q(2)$, where we
can apply Theorem \ref{psuq}. As it is the case for the other computations of Poisson boundaries of quantum groups in the literature, this is the
most subtle part of the proof.
\item Once the multiplicativity of $\iz$ on $\Sone \frakB$ is proved, the general results gathered in Proposition \ref{prop.abelian} allow to
conclude.
\end{itemize}

\subsection*{Notations and equivariance formulae}

\begin{notation} \label{not.intertwiners}
All objects related to $\cGtil=\SU_q(2)$ are denoted with tildes, while the corresponding objects related to $\cG = A_o(F)$ are denoted without
tildes. Fix a monoidal equivalence $\vphi : \cG \recht \cGtil$. We identify $\Irred(\cG) = \Irred(\cGtil) = \N$ and we make once and for all a choice
of \emph{isometric} intertwiners $\Vtil(x \ot y,z) \in \Mor(x \ot y,z)$ for $\cGtil=\SU_q(2)$. We take $V(x \ot y,z)$ such that $\vphi(V(x \ot y,z))
= \Vtil(x \ot y,z)$.

Whenever $z \in x \ot y$, we set $p^{x \ot y}_z = V(x \ot y,z) V(x \ot y,z)^*$ and we define $\ptil^{x \ot y}_z$ accordingly. For any $x \in \N$, we
have irreducible representations $U^x$ of $\cG$ on $H_x$ and $\Util^x$ of $\cGtil$ on $\Htil_x$. Finally, recall from \ref{not.states}, the special
states $\psi_x$ on $\B(H_x)$ (and hence, $\psitil_x$ on $\B(\Htil_x)$).
\end{notation}

\begin{proposition} \label{prop.valuesinpoisson}
Let $\iz$ be as in \eqref{eq.geniz} and let $\mu$ be a probability measure on $\Irred(\cG)$. Then, $\iz(a) \in \hon$ for all $a \in \frakB$.
\end{proposition}
\begin{proof} Let $a\in {\mathfrak B}$. Then,
for $x,y\in \mbox{Irred}(\mathbb{G})$, we know that
$$(p_x\ot p_y)\deh(\iz(a))=\sum_{z\in x\ot y}V(x\ot y, z)(\iz(a)p_z)V(x\ot y, z)^*$$
So
\begin{align}
(p_x\ot p_y)\deh(\iz(a)) &= \sum_{z\in x\ot
  y}(\id\ot\id\ot\omega)\bigl((V(x\ot y, z)\ot 1)(X^z)^*(1\ot a)X^z
(V(x\ot y, z)^*\ot 1)\bigr) \notag \\ &=
\sum_{z\in x\ot y}(\id\ot\id\ot\omega)\bigl((X_{23}^y)^*(X_{13}^x)^*(\widetilde{p}_{z}^{x\ot y}\ot a
)X_{13}^xX_{23}^y \bigr)\label{iz1}\\
&= (\id\ot\id\ot\omega)\bigl((X_{23}^y)^*(X_{13}^x)^*(1\ot 1\ot a )X_{13}^xX_{23}^y \bigr)\notag
\end{align}
where \eqref{iz1} is valid because $(V(x\ot y, z)\ot
1)(X^z)^*=(X^y)_{23}^*(X^x)_{13}^*(\widetilde{V}(x\ot y, z)\ot 1)$. Then, we get that
\begin{align*}
p_x(\id\ot \psi_y)\deh(\iz(a)) &
= (\id\ot \psi_y\ot \omega)\bigl((X_{23}^y)^*(X_{13}^x)^*(1\ot 1\ot a )X_{13}^xX_{23}^y  \bigr) \\
&= (\id\ot \omega)\bigl((X^x)^*(1\ot a)X^x\bigr) = \iz(a)p_x
\end{align*}

The equality follows from the fact that $(\psi_y\ot \omega)\bigl((X^y)^*(1\ot b)X^y \bigr)=\omega(b)$ for $b\in \wB$. Indeed, the KMS-property of the
state $\omega$ (see \eqref{KMSS}) gives
\begin{equation*}
(\psi_y\ot \omega)\bigl((X^y)^*(1\ot b)X^y \bigr) = (\widetilde{\psi}_y\ot \om)\bigl((\widetilde{Q}_y^{-2}\ot 1)X^y (X^y)^*(1\ot b)\bigr) =
(\widetilde{\varphi}_y\ot \om)(1\ot b)  = \om (b)
\end{equation*}
This completes the proof.
\end{proof}

\begin{proposition} \label{prop.intertwines}
The generalized Izumi operator $\iz$ defined by \eqref{eq.geniz} intertwines the actions $\beta_{\cG},\beta_{\cGh}$ (see Def.\ \ref{def.B} and
formula \eqref{eq.adjoint}) with the actions $\al_{\cG},\al_{\cGh}$ defined in \ref{def.actions}.
\end{proposition}
\begin{proof}
Intertwining of $\al_\cG$ and $\be_\cG$ follows from
\begin{align*}
(\iz\ot\id)\beta_{\cG}(a) &= (\id\ot \omega\ot\id)(\cX_{12}^*(1\ot \beta_{\cG}(a))\cX_{12})
= (\id\ot \omega\ot\id)\bigl(\mathbb{V}_{13}(\id\ot\beta_{\cG})\bigl(\cX^*(1\ot a)\cX \bigr)\mathbb{V}_{13}^*\bigr) \\
&= \mathbb{V}\bigl((\id\ot \omega\ot\id)(\id\ot \beta_{\cG}) \bigl(\cX^*(1\ot a)\cX \bigr)\bigr)\mathbb{V}^*
= \mathbb{V}\bigl((\id\ot \omega) \bigl(\cX^*(1\ot a)\cX \bigr)\ot 1\bigr)\mathbb{V}^* \\
&= \al_\cG(\iz(a)) \; .
\end{align*}
On the other hand, intertwining of $\al_{\cGh}$ and $\be_{\cGh}$ is a consequence of
\begin{align*}
\al_{\cGh}(a) V(y \ot z,x) &= \widehat{\de}(\iz(a))V(y\ot z,x) = V(y\ot z,x)\iz(a)p_x
= V(y\ot z,x)(\id \ot \om)((X^x)^*(1\ot a)X^x) \\
&= (\id\ot\id\ot\om)\bigl((X^z)^*_{23}(X^y)^*_{13}(1\ot 1\ot a)X^y_{13}X^z_{23}\bigr)V(y\ot z,x) \\
&= (\id\ot \iz)(\cX^*(1\ot a)\cX)V(y\ot z,x) = (\id\ot \iz)\be_{\cGh}(a)V(y\ot z,x)
\end{align*}
for every $x,y,z\in \Irred(\cG)$.
\end{proof}

\subsection*{Multiplicativity of the generalized Izumi operator}

Fix a generating probability measure $\mu$ on $\Irred(\cG)$. We prove that the generalized Izumi operator $\iz$ given by \eqref{eq.geniz}, is
multiplicative on $\Sone \frakB$, using a strategy from \cite{INT}.

For every $y\in \mbox{Irred}(\mathbb{G})$, we define
$$\iz_y:\mathfrak{B}\recht B(H_y):\iz_y(a):=\iz(a) p_y=(\id\ot\om)((X^y)^*(1\ot a)X^y)$$

\begin{lemma} \label{lem.adjoint}
For the scalar products on $\mathfrak{B}$ and $B(H_y)$ given respectively by $\om$ and $\psi_y$, the linear map $\iz_y$ has an  adjoint mapping
$\iz_y^*:B(H_y)\recht \mathfrak{B}$ given by
$$\iz_y^*(b)=(\widetilde{\varphi}_y\ot  \id)(X^y(b\ot 1)(X^y)^*) $$
So $\om(a\iz_y^*(b))=\psi_y(\iz_y(a)b)$ for $a\in \mathfrak{B}$ and $b\in B(H_y)$.
\end{lemma}
\begin{proof} Again, we use the KMS-property of $\om$. From this it follows that
\begin{align*}
\psi_y(\iz_y(a)b) &= \psi_y\bigl(\id\ot \om) ((X^y)^*(1\ot a)X^y b) \bigr) = (\psi_y\ot\om)((X^y)^*(1\ot a)X^y (b\ot 1))\\
&= (\widetilde{\varphi}_y\ot \om)((1\ot a)X^y(b\ot 1)(X^y)^*) = \om\bigl(a(\widetilde{\varphi}_y\ot  \id)(X^y(b\ot 1)(X^y)^*) \bigr)\\
&= \om(a\iz_y^*(b))
\end{align*}
\end{proof}

\begin{remark}
We already know that $(\iz_y\ot \id)\be_\cG=\al_\cG\circ\iz_y$. From the definition of the adjoint $\iz_y^*$ it follows that $(\iz_y^*\ot
\id)\al_\cG=\be_\cG\circ \iz_y^*$. This follows also from the next calculation.
\begin{align*}
(\iz_y^*\ot \id)\al_\cG(a) &= (\widetilde{\varphi}_y\ot\id \ot\id)\bigl(X_{12}^y(\al_\cG(a))_{13} (X_{12}^y)^*\bigr) = (\widetilde{\varphi}_y\ot\id \ot\id)(X_{12}^y V_{13}(a\ot 1\ot 1)V_{13}^*(X_{12}^y)^*) \\
&= (\widetilde{\varphi}_y\ot\id \ot\id)\bigl((\id\ot \be_\cG)(X^y(a\ot 1)(X^y)^*)\bigr)
= \be_\cG\bigl((\widetilde{\varphi}_y\ot\id)(X^y(a\ot 1)(X^y)^*) \bigr) \\
&= \be_\cG\bigl(\iz_y^*(a) \bigr)
\end{align*}
\end{remark}

\begin{lemma} \label{lem.in-mult}
Let $\mu$ be a generating probability measure on $\Irred(\cG)$. Denote
$$\rP_y=\iz_y^*\iz_y:\mathfrak{B}\recht \mathfrak{B}\qquad
\mbox{and}\qquad \rP_{\mu}=\sum_{y}\mu(y)\rP_y \; .$$
The sequence $(\rP_{\mu^{* n}})_n$ converges pointwize $^*$-strongly
to a completely positive unital map $\rP_\infty : \frakB \recht
\frakB$. Moreover, if $a \in \frakB$, the following are equivalent.
\begin{itemize}
\item The element $a$ belongs to the multiplicative domain of the completely positive unital map $\iz$.
\item We have $\rP_\infty(a) = a$.
\end{itemize}
\end{lemma}
\begin{proof}
We first make the following claim.

{\bf Claim.} For every $a \in \frakB$, the sequence $\rP_{\mu^{*n}}(a)$ converges $^*$-strongly. Denoting its limit as $\rP_\infty(a)$, we have
$$\epsh(\iz(b) \cdot \iz(a)) = \om(b \rP_\infty(a)) \quad\text{for all}\;\; a,b \in \frakB \; .$$
We prove the claim below and argue already how the lemma follows from it. If $a \in \frakB$ belongs to the multiplicative domain of $\iz$, we have
$\iz(b) \cdot \iz(a) = \iz(ba)$ for all $b \in \frakB$. The formula in the claim yields $\om(ba) = \epsh(\iz(ba))=\om(b \rP_\infty(a))$ for all $b
\in \frakB$. Hence, $\rP_\infty(a) =a$. Conversely, if $\rP_\infty(a) = a$, we find that $\epsh(\iz(a)^* \cdot \iz(a)) = \epsh(\iz(a^*a))$. Since
$\iz(a)^* \cdot \iz(a) \leq \iz(a^* a)$ and since $\epsh$ is faithful on $H^\infty(\cGh,\mu)$, it follows that $\iz(a)^* \cdot \iz(a) = \iz(a^*a)$.
This implies that $a$ belongs to the multiplicative domain of $\iz$.

It remains to prove the claim. For all $a,b \in \frakB$, we have
\begin{align*}
\widehat{\varepsilon}(\iz(b)\cdot\iz(a)) &= \widehat{\varepsilon}\bigl(\lim_{n\recht \infty}P_{\mu^{* n}}(\iz(b)\iz(a)) \bigr)
=  \lim_{n\recht \infty} \psi_{\mu^{* n}}(\iz(b)\iz(a)) \\
&= \lim_{n\recht \infty}\sum_{x}\mu^{* n}(x)\psi_x(\iz_x(b)\iz_x(a)) = \lim_{n\recht \infty}\om\bigl(b\rP_{\mu^{* n}}(a)\bigr) \; .
\end{align*}
It follows that $\rP_{\mu^{*n}}(a)$ is weakly convergent, say to $\rP_\infty(a) \in \frakB$. Since all the $\rP_{\mu^{*n}}$ and $\rP_\infty$ commute
with the ergodic action $\be_\cG$ of $\cG$ on $\frakB$ and preserve the state $\om$, these completely positive operators preserve the
(finite-dimensional) spectral subspaces of $\be_\cG$ and it follows that $\rP_{\mu^{*n}}(a) \recht \rP_\infty(a)$ $^*$-strongly. This proves the
claim.
\end{proof}

With respect to the ergodic action $\be_\cG$, the von Neumann algebra $\frakB$ has the natural dense $^*$-subalgebra $\cB$ given as the linear span
of the spectral subspaces $B_x$, $x \in \Irred(\cG)$. Since $\cP_y$ commutes with $\be_\cG$, it follows that $\cP_y$ maps $B_x$ into $B_x$. We study
more closely this operator on the finite-dimensional vector space $B_x$.

Denote
$$\iz_y^x:B_x\recht B(H_y) : \iz_y^x(a) = \iz(a)p_y$$
the restriction of the operator $\iz_y$ to the
spectral subspace $B_x$. We use the following unitary identifications
that are consequences of Notation \ref{not.states} and formula \eqref{eq.orthogonalitymon}.
\begin{align}
v_y & :(B(H_y),\psi_y)\recht H_y\ot H_y: A \mapsto
\frac{1}{\sqrt{\dimq(y)}}(A\ot 1)t_y \notag \\
\phi_x & : (B_x,\om)\recht H_x\ot \overline{\Htil}_x:
  (\om_{\widetilde{\mu},\rho}\ot\id)(X^x) \mapsto
  \frac{1}{\sqrt{\dimq(x)}}\ \rho\ot
  \overline{\widetilde{Q}^{-1/2}_x\widetilde{\mu}} \label{eq.ident}
 \end{align}
with $\rho\in H_x$ and $\widetilde{\mu}\in \Htil_x$.

\begin{lemma} \label{lem.computation}
For every $x,y \in \Irred(\cG)$, we have
$$v_y\circ\iz_y^x\circ \phi_x^*=
V(y\ot y,x)(1 \ot {\overline{\xi_y^x}}^*)
$$
where the vector $\xi_y^x \in \Htil_x$ is defined as
$$
{{\xi}_y^x} =\frac{1}{\sqrt{\dimq(x)\dimq(y)}}
 \Vtil(y\ot y,x)^* (\Qtil_y^{-2}\ot 1)\ttil_y \; .
$$
\end{lemma}
\begin{proof}
Take $a\in B(H_y)$ and $b\in B_x$. By Lemma \ref{lem.adjoint}, we have
$$
\psi_y\bigl(a^*\iz_y^x(b)\bigr) = \om\bigl((\widetilde{\varphi}_y\ot
\id)(X^y(a\ot 1)(X^y)^*)^*b\bigr) \; .
$$
Further calculation gives
\begin{align*}
(\widetilde{\varphi}_y\ot \id)\bigl(X^y(a\ot 1)(X^y)^*\bigr)
&= (\widetilde{\psi}_y\ot \id)\bigl( (\Qtil_y^{-2}\ot 1) X^y (a\ot 1)
(X^y)^* \bigr) \\
&= \frac{1}{\dimq(y)}(\ttil_y^*\ot 1) (\Qtil_y^{-2}\ot 1\ot 1)
X^y_{13}(a\ot 1\ot 1)(X^y_{13})^* (\ttil_y\ot 1) \\
&= \frac{1}{\sqrt{\dimq(y)}} (\ttil_y^*\ot 1) (\Qtil_y^{-2}\ot 1\ot 1)
X^y_{13}X^y_{23} (v_y(a) \ot 1) \; .
\end{align*}
We know that
\begin{equation}\label{eq.b}
X^y_{13}X^y_{23}=\sum_{z\in y\ot y}\bigl(\Vtil(y\ot y,z)\ot
1\bigr)X^z\bigl({V}(y\ot y,z)^* \ot 1\bigr) \; .
\end{equation}
So, $\psi_y\bigl(a^*\iz_y^x(b)\bigr)=\om (D^* b)$ where
\begin{equation*}
D=\frac{1}{\sqrt{\dimq(y)}}(\ttil_y^* (\Qtil_y^{-2} \ot 1) \Vtil(y \ot
y,x) \ot 1) \; X^x \; (V(y \ot y,x)^* v_y(a) \ot 1) \; .
\end{equation*}
Here, only the term $z=x$ in the sum \eqref{eq.b} remained because of the orthogonality relations \eqref{eq.orthogonalitymon} and the assumption
$b\in B_x$.

Using the commutation relations $(\Qtil_y^{-1/2}\ot \Qtil_y^{-1/2})\Vtil(y\ot y,x)=\Vtil(y\ot y,x) \Qtil_x^{-1/2}$
and $(\Qtil_y^{-1/2}\ot \Qtil_y^{-1/2}) \ttil_y=\ttil_y$ and the
formula $\om(D^* b) = \la \phi_x(D),\phi_x(b) \ra$, we arrive at
$$\psi_y(a^* \iz_y^x(b)) = \la V(y \ot y,x)^* v_y(a) \ot
\overline{\xi^x_y} , \phi_x(b) \ra$$
where $\xi^x_y$ is given in the statement of the lemma. This proves
the lemma.
\end{proof}

We finally prove the multiplicativity of the operator $\iz$ on $\Sone \frakB$.

\begin{lemma} \label{lem.multiplicative}
The elements of $\Sone \frakB$ belong to the multiplicative domain of the generalized Izumi operator $\iz$ introduced in \eqref{eq.geniz}.
\end{lemma}
\begin{proof}
Since the action $\be_\cG$ on $\frakB$ commutes with the action $\rho$
of $S^1$ on $\frakB$, it suffices to show that every element of $B_x$
that is invariant under $\rho$ belongs to the multiplicative domain of
$\iz$.

Observe that for $z = |q|^{it}$, we have $v_x \rho_z = (1 \ot \overline{\Qtil^{-it}}) v_x$. Combining Lemmas \ref{lem.in-mult} and \ref{lem.computation}, the
multiplicativity of $\iz$ on $\Sone \frakB$ is then equivalent with the statement
$$\lim_{n \recht \infty} \sum_y \mu^{*n}(y) \; \xi_y^x \; \la \xi_y^x, \eta \ra \recht \eta \quad\text{for all}\;\; \eta \in \Htil_x^{S^1} \; ,$$
where we denote by $\Htil_x^{S^1}$ the subspace of vectors $\eta \in \Htil_x$ satisfying $Q_x \eta = \eta$.

This last statement concerns only $\SU_q(2)$ and hence holds because of Theorem \ref{psuq}.
\end{proof}

\subsection*{Proof of Theorem \ref{thm.main}}

We have gathered enough material to give the

\begin{proof}[Proof of Theorem \ref{thm.main}]
By Proposition \ref{prop.valuesinpoisson} and Lemma \ref{lem.multiplicative}, we have the normal unital $^*$-homomorphism
$$\iz : \Sone \frakB \recht \hon \; .$$
Because $\om$ is faithful, the map $\iz$ is injective and by Proposition \ref{prop.intertwines}, $\iz$ intertwines the actions
$\be_{\cG}$, $\be_{\cGh}$ on $\Sone \frakB$ with the actions $\al_\cG$, $\al_{\cGh}$ on $\hon$. So, it remains to prove the surjectivity of $\iz$.

From Proposition \ref{prop.abelian}, we know that $\al_\cG$ is an ergodic action of $\cG$ on $\hon$ and that $\mult(x,\al_\cG) = 0$ for $x$ odd,
$\mult(x,\al_\cG) \leq 1$ for $x$ even. On the other hand, the multiplicitity of $x$ in the restriction of $\be_\cG$ to $\Sone \frakB$ is exactly $1$
for $x$ even. Indeed, as we have seen in the proof of Lemma \ref{lem.multiplicative}, $v_x$ maps the $S^1$-invariant elements of $B_x$ onto $H_x \ot
\overline{\Htil_x^{S^1}}$ and $\Htil_x^{S^1}$ is one-dimensional for
$x$ even. Indeed, the eigenvalues of $\Qtil_x$ are
$q^x,q^{x-2},\ldots,q^{2-x},q^{-x}$, all of multiplicity $1$ on the Hilbert
space $\Htil_x$ of dimension $x+1$.
Injectivity and equivariance of $\iz$ finish the proof.
\end{proof}

\section{Martin boundary of the dual of $\aof$}

We prove an identification theorem for the Martin boundary
of the dual of $\aof$.

As above, fix $\cG = \aof$ with $\cG \neq
\SU_{\pm 1}(2)$. Take $q$, the \cst-algebra $B$ and the actions $\be_\cG$ and
$\be_{\cGh}$ as in Definition \ref{def.B} and formula
\eqref{eq.adjoint}. The \cst-algebra $B$ admits the natural action
$(\rho_z)_{z \in S^1}$ of $S^1$ and the subalgebra of $S^1$-invariant
elements was denoted by $\Sone B$.

Recall from \eqref{eq.actionsMartin} that the
Martin boundary $A_\mu$ is naturally equipped with actions
$\gamma_\cG$ and $\gamma_{\cGh}$.

\begin{theorem} \label{thm.martin}
Let $\mu$ be a generating measure on $\Irred(\cG)$ that is transient
and has finite first moment:
$$\sum_{x \in \N} x \mu(x) < \infty \; .$$
The composition of the generalized Izumi operator given in
\eqref{eq.geniz} and the quotient map $\pi : \ell^\infty(\cGh) \recht
\frac{\ell^\infty(\cGh)}{c_0(\cGh)}$ defines a $^*$-isomorphism
$$\pi \circ \iz : \Sone B \recht A_\mu$$
of $\Sone B$ onto the Martin boundary $A_\mu$. This $^*$-isomorphism
intertwines the actions $\be_{\cG},\be_{\cGh}$ on $\Sone B$ with the actions
$\gamma_\cG,\gamma_{\cGh}$ on $A_\mu$.
\end{theorem}
\begin{proof}
Combining Theorems 5.6 and 5.8 in \cite{Vbo}, $A_\mu$ admits a state
$\om_\infty$ such that the map
$$T : A_\mu \recht \hon : T(a) = (\id \ot
\om_\infty)\gamma_{\cGh}(a)$$
is a $^*$-homomorphism with dense range. Moreover, by formula (5.3) in
\cite{Vbo}, the composition $\pi \circ T$ is the identity map. Since
$T$ intertwines the action $\gamma_\cG$ on $A_\mu$ with the action
$\al_\cG$ on $\hon$, it follows that $T$ is a $^*$-isomorphism of
$A_\mu$ onto the \cst-algebra $D$ defined as the closed linear span of the spectral subspaces of the action
$\al_\cG$ of $\cG$ on $\hon$. Moreover, the quotient map $\pi$
provides the inverse of this $^*$-isomorphism.

By Theorem \ref{thm.main}, the restriction of the generalized Izumi
operator, yields the $^*$-isomorphism $\iz : \Sone B \recht
D$. Composing $\pi$ and $\iz$, we are done.
\end{proof}

\end{document}